\documentclass[12pt]{amsart}

\newtheorem{Theorem}{Theorem}[section]
\newtheorem{Proposition}[Theorem]{Proposition}
\newtheorem{Def}[Theorem]{Definition}
\newtheorem{Lemma}[Theorem]{Lemma}\newtheorem{Example}[Theorem]{Example}
\newtheorem{Corollary}[Theorem]{Corollary}

\newcommand{\el}{\mbox{$\mathcal L$}}
\newcommand{\dee}{\mbox{$\mathcal D$}}
\newcommand{\sqs}{\mbox{${\mathcal S}(S)$}}

\newcommand{\jay}{\mbox{$\mathcal J$}}
\newcommand{\pee}{\mbox{$\mathcal P$}}
\newcommand{\ar}{\mbox{$\mathcal R$}}
\newcommand{\eh}{\mbox{$\mathcal H$}}
\newcommand{\leqel}{\mbox{$\leq _{\mathcal L }$}}
\newcommand{\leqjay}{\mbox{$\leq _{\mathcal J }$}}
\newcommand{\leqar}{\mbox{$\leq _{\mathcal R}$}}
\newcommand{\leqels}{\mbox{$\leq _{{\mathcal L}^{\ast}}$}}
\newcommand{\leqars}{\mbox{$\leq _{{\mathcal R}^{\ast}}$}}

\newcommand{\els}{\mbox{${\mathcal L}^{\ast}$}}

\newcommand{\ars}{\mbox{${\mathcal R}^{\ast}$}}
\newcommand{\ehs}{\mbox{${\mathcal H}^{\ast}$}}

\newcommand{\s}{\mbox{$^{\sharp}$}}

\begin{document}

\title[Semigroups of quotients]
{Semigroups of left quotients:\\ existence, straightness and locality}
\subjclass{20 M 10}

\keywords{Group inverse, semigroup of (left) quotients,
order, straightness, locality}

\author{Victoria Gould}
\address{Department of Mathematics\\University of York\\Heslington\\York
YO10 5DD\\UK}
\email{varg1@york.ac.uk}
\date{\today}

\begin{abstract}

A subsemigroup $S$ of a semigroup $Q$ is a {\em local left
order} in $Q$ if, for every group $\eh$-class $H$ of $Q$,
$S\cap H$ is a left order in $H$ in the sense of group theory.
That is, every $q\in H$ can be written as $a\s b$ for some $a,b\in S\cap H$,
where $a\s$ denotes the group inverse of $a$ in $H$.
On the other hand,  $S$ is a  {\em left order} in $Q$
and $Q$ is a {\em semigroup of left quotients} of $S$
 if every element of $Q$
can be written as $c\s d$
where $c,d\in S$ and if,
in addition, every element of $S$ that is
{\em square cancellable} lies in a subgroup
of $Q$.   If one also
insists  that $c$ and $d$ can be chosen such that
$c\,\ar\, d$ in $Q$, then $S$ is said to be a {\em straight left order}
 in $Q$.

This paper investigates the close relation between local left
orders and straight left orders in a semigroup $Q$
and gives some quite general conditions for a left
order $S$ in $Q$ to be straight. In the light of the
connection between locality and straightness
we give a complete description of straight left
orders that improves upon that in our earlier paper.
\end{abstract}
\maketitle
\section{Introduction}

The concepts of left order in ring theory and in group theory have
proved fruitful in the investigation of non-regular rings
and cancellative semigroups.  This paper is concerned with the analogous
notion
for semigroups, introduced by Fountain and Petrich in
\cite{fpi}. Their aim was to develop concepts that reflect the equal
importance of {\em all} subgroups of a semigroup, not only the group
of units, which of course may not even exist.
A subsemigroup $S$ of a semigroup $Q$ is a {\em left order} in $Q$
and $Q$ is a {\em semigroup of left quotients} of $S$ if every element
of $Q$ can be written as $a\s b$ for some $a,b\in S$ and if, in
addition,
every element of $S$ that satisfies a weak cancellation property
called {\em square cancellation}, lies in a subgroup 
of $Q$. Here
$a\s$ denotes the group inverse of $a$ in a subgroup of $Q$. {\em Right orders } and {\em semigroups
of right quotients} are defined dually; if $S$ is both a left order
and a right order in $Q$ then $S$ is an {\em  order} in $Q$ and
$Q$ is a {\em semigroup of  quotients} of $S$.

In the case where $S$ is a subring of  a ring $Q$ our concepts are
closely related to the classical ones of ring theory \cite{rings}. 
Ore's theorem states that a ring $R$ has a ring of left quotients if
and only if $R$ contains a non-zero divisor  and for all
$s,t\in R$ where $s$ is a non-zero divisor, there exist $h,k\in R$ with
$h$ a non-zero divisor and $ht=ks$ \cite{l}.

In the case where $S$ is a subsemigroup of
a group $Q$, our concept of a semigroup of left
quotients becomes the classical one of a group of left quotients.
A theorem of Ore and Dubreil \cite{cp} says
 that  a semigroup $S$ has a group of left quotients
if and only if it is right reversible, that is, for all $a,b\in S$,
$Sa\cap Sb\neq\emptyset$, and cancellative. 

In view of the above, the question then arises
of characterising left orders in semigroups from
larger classes than the class of groups. Of course, one chooses classes
whose semigroups have a structure that is well understood, the point
being that the left orders inherit some of this structure. 
This path is well trodden: for a sample we
 refer the reader to \cite{fpi} and \cite{compreg}.
Results of this kind 
make it  natural to
 ask the more
searching question  of whether we can find an analogue for semigroups
of the theorems of Ore and Dubreil.  In its most general
form this would appear to be very difficult, the reason being it is
hard to have control over products of quotients $(a\s b)(c\s d)$.
However, in many cases, it is always possible to write quotients in the
form $a\s b$ where $a\,\ar\, b$ in the semigroup of left quotients under
consideration. For example, this is the case if the latter is a regular
semigroup on which Green's relation $\eh$ is a congruence. If we can
write quotients in this way we 
have more idea where (in terms of the $\dee$-class structure
of $Q$) products of quotients will lie.
This brings us to {\em straightness}.

We say that a left order $S$ in $Q$ is {\em straight} if every element
of $Q$ can be written as $a\s b$ for some $a,b\in S$ with $a\,\ar\, b$
in $Q$. If $S$ is straight in $Q$ it is easy to see that
$Q$ is regular and $S$ is {\em very large}
in $Q$, which means that $S$ has non-empty intersection with every
$\eh$-class of $Q$.  Consequently, $S$ inherits much of the
$\dee$-class structure of $Q$. Moreover $S$ is {\em local}
in $Q$ in the sense that each group $\eh$-class $H$ of $Q$
is the group of left quotients of $H\cap S$.  By Ore and Dubreil's
theorem, the subsemigroups $S\cap H$ are right reversible and cancellative.

Straight left orders were first investigated in \cite{vienna}, where
the above properties are proven.  The main theorem of \cite{straight}
gives a complete description of straight left orders. This
is achieved via 
{\em embeddable $*$-pairs}, which are pairs of preorders on a semigroup
$S$ that have properties reflecting those of $\leqel$ and $\leqar$ on
a regular semigroup. To be a little more specific, Theorem 4.1 of
\cite{straight} states that if $(\leq_l,\leq_r)$
is an embeddable $*$-pair for a semigroup $S$, then
$S$ is a straight left order in a  (regular) semigroup
$Q$  such that for all $b,c\in S$,
\[b\leq_l c\mbox{ in }S\mbox{ if and only if }b\leqel c\mbox{ in } Q\]
and
\[b\leq_r c\mbox{ in }S\mbox{ if and only if }b\leqar c\mbox{ in } Q\]
if and  only if $S$ satisfies four given conditions.
The first, (Gi), ensures that square cancellable elements of
$S$ lie in subgroups of $Q$ and the second, (Gii),
states that certain subsemigroups of $S$ are right reversible;
thus far the conditions look natural. The remaining
two conditions, (Giii) and (Giv),
 are necessary for the construction of $Q$ from $S$
but are aesthetically unpleasing. For the details we refer the reader
to \cite{straight} and Section 4 below. The main task of
this paper is to show that conditions (Giii) and (Giv) {\em always}
hold for an embeddable $*$-pair satisfying (Gi) and (Gii). This
yields a much cleaner description of straight left orders that
is moreover far easier to use, as we demonstrate in Section 5
and more significantly in \cite{appl} and \cite{poorders}.

The layout of the paper is as follows.
In Section 2 we recall some details of the preorders
$\leqel$ and $\leqar$ and the related preorders
$\leqels$ and $\leqars$. We also discuss square cancellability
and give the formal definitions of left order and variations
of such.

In Section 3 we consider the relationship between local left
orders and straight left orders. We give some examples of orders
that are not straight. On the other hand we find some quite
general conditions
which ensure straightness. 

Section 4 contains the promised streamlined description of
straight left orders. For convenience we also give the details
of embeddable $*$-pairs. We pick our way carefully through related
properties with an eye to applications, to appear in \cite{appl}.

The final section gives some straightforward specialisations of
Theorem 4.4. Specifically, we 
characterise straight left orders in inverse semigroups, in
completely regular semigroups and more generally in completely
semisimple semigroups.

\section{Preliminaries}

 This section recalls properties of
certain preorders on a semigroup $S$. The behaviour
of these preorders is 
 essential in determining the existence and in working with
semigroups of (left) quotients of $S$.

We assume a familiarity with the basic notions of semigroup theory, as
in
\cite{cp} and  \cite{h}, in
particular, with Green's relations. We bear in mind the fact that
Green's relations $\el$ and $\ar$ are the equivalence relations
associated with the  preorders $\leqel$ and $\leqar$ respectively,
where for elements $p$ and $q$ of a semigroup $Q$,
 \[p\,\leqel\, q\mbox{ if and only if }Q^1p\subseteq Q^1q\]
and dually,
 \[p\,\leqar\, q\mbox{ if and only if }pQ^1\subseteq qQ^1.\]
Notice that the relations
 $\leqel$ and $\leqar$ are respectively  right and
left compatible  with multiplication.

For any element $q$ of a semigroup $Q$, $q\s$ denotes the inverse of
$q$ in a subgroup of $Q$. Use of the notation $q\s$ implies that
$q$ lies in a subgroup of $Q$; $q\s$ is then uniquely determined as the
(group) inverse of $q$ in the subgroup $H_q$ of $Q$. The union of
subgroups of $Q$ is denoted by $\mathcal{H}(Q)$.

For convenience we list some properties of $\leqel$
and $\leqar$ that will be used repeatedly. 

\begin{Lemma} Let $p,q,r$ be elements of a semigroup $Q$. The following
statements and their left-right duals hold:

(i) if $q\,\leqel\, p$ and $p\in\eh (Q)$, then $qp\,\ar\, q$;

(ii) if $q,r\,\leqel\, p, p\in\eh(Q)$ and $qp=rp$, then $q=r$;

(iii) if $q,r\, \leqel\, p,p\in\eh(Q)$ and $qp\,\el\, rp$, then $q\,\el\, r$;

(iv) if $q\,\leqar\, p$ and $q\in\eh (Q)$, then $qp\,\ar\, q$;

(v) if $q\,\ar\, p$ and $p\in\eh (Q)$, then $pq\,\eh q$.
\end{Lemma}
\begin{proof} Conditions (i) to (iv)
are taken from Lemma 2.1 of \cite{compreg}. For (v)
notice that the dual of (i) gives that $pq\,\el\, q$.
Since $\ar$ is a left congruence we have
\[pq\,\ar\, p^2\,\eh\, p\,\ar\, q\]
so that $pq\,\eh\, q$ as required.
\end{proof}

We make repeated use of two further preorders, namely
$\leqels$ and $\leqars$. For elements $a$ and $b$ of a semigroup $S$,
$a\,\leqels\, b$ if and only if for all $x,y\in S^1$,
\[bx=by\mbox{ implies that }ax=ay.\]
Clearly $\leqels$ is a preorder that is right
compatible with multiplication so that
the associated equivalence relation,
denoted $\els$, is a right congruence. The relations
$\leqars$ and $\ars$ are defined dually.

Our next lemma is well known and can easily be proved by embedding
$S$ in the (dual of the)
full transformation semigroup on the set $S^1$.

\begin{Lemma} \cite{f},\cite{L} Let $S$ be a semigroup and let $a,b\in S$. Then
$a\,\leqels \,b$ if and only if $a\,\leqel \,b$ in an oversemigroup of $S$.
The corresponding statements are true for the relations $\leqars,
\els$ and $\ars$.
\end{Lemma}

By Green's Theorem (Theorem 2.2.5 of \cite{h}),
$q$ lies in a subgroup (so that $q\s$ exists)
if and only if $q\,\mathcal{H}\, q^2$. Denoting by $\ehs$ the
intersection
of $\els$ and $\ars$, it is then clear that a necessary condition 
for an element $a$ of a semigroup $S$ to lie in a subgroup of
an oversemigroup is that $a\,\ehs\, a^2$.
We say that an element $a$ of a semigroup $S$ is {\em square
cancellable}
if $a\,\ehs\, a^2$ and we denote by $\sqs$ the set of square cancellable
elements of $S$.
We can now give the definitions of our core concepts.

\begin{Def} A subsemigroup $S$ of a semigroup $Q$ is
a {\em  left order } in $Q$  if every element of $Q$ can
be written as $a\s b$ for some $a,b\in S$
and if, in addition,
every $a\in\sqs$ lies in a subgroup of $Q$.  If $S$ is a left
order in $Q$, then we also say that
$Q$ is a {\em semigroup of left quotients} of $S$.

The concepts of {\em right order} and {\em semigroup
of right quotients} are the appropriate duals. If $S$
is both a left order and a right order in $Q$
then we say that $S$ is an {\em order} in
$Q$ and $Q$ is a {\em semigroup of quotients} of $S$.
\end{Def}

In the above definitions, if we omit the condition
that square cancellable elements must lie in subgroups, then
we insert the adjective `weak'. 

It is easy to see that if $S$ is a weak left order in $Q$, then any
element of $Q$ may be written as $a\s b$ where $a,b\in S$ and
$b\,\leqar\, a$ in $Q$. In this case the dual of
Lemma 2.1 (i) gives that
$a\s b\,\el\, b$ in $Q$. Thus $S$ intersects
 every $\el$-class of $Q$.  If moreover $S$ intersects
every $\eh$-class of $Q$, then $S$ inherits more of the structure of
$Q$. One way of ensuring this is by utilising the
 concept of {\em straightness}.

Let $S$ be a (weak) left order in $Q$. Then $S$ is {\em straight }
in $Q$ if every element of $Q$ can be written as $a\s b$ where
$a,b\in S$ and $a\,\ar\, b$ in $Q$. For such elements
$a$ and $b$ Lemma 2.1 (v) gives that $a\s b\,\eh\, b$ in $Q$, so
that $S$ intersects every $\eh$-class of $Q$. Note also that
$Q$ must be regular. 
The adjective `straight' can be inserted in any 
of the concepts in Definition 2.3. However, some
care should be taken: a right order $S$ in $Q$ is
straight if every element of $Q$ can be written as
$ba\s$ for some $a,b\in S$ with $a\,\large{\mathcal{L}}\, b$ in $Q$.
A straight order is an order that is straight as a left
order and as a right order.

Straightness is a curious property, naturally occuring for many
left orders. For example, if 
$Q$ is a regular semigroup for which Green's relation $\eh$ is
a congruence, then every left order in $Q$ is straight \cite{vienna}.
On the other hand a completely regular semigroup can contain
left orders which are {\em not} straight; however, every 
order in a completely regular semigroup {\em is} straight
\cite{compreg}. Straightness is also related to locality.

\begin{Def} Let $S$ be a subsemigroup of a semigroup $Q$. Then
$S$ is a {\em local left order} in $Q$ if, for every group
$\eh$-class $H$ of $Q$, $S\cap H$ is a left order in $H$.
{\em Local right orders} are defined dually.
\end{Def}

The
next section attempts to throw some light on these connections.

\section{Straightness, locality and size}
This section investigates the relationship between local
left orders, straight left orders and {\em very large} left orders.
By a very large subsemigroup of a semigroup $Q$ we mean one that has
non-empty 
intersection with every $\eh$-class of $Q$.

\begin{Proposition}
Let $S$ be a subsemigroup of a regular semigroup $Q$. The
following conditions are equivalent:

(i) $S$ is a very large weak left order in $Q$:

(ii) $S$ is a weak straight left order in $Q$;

(iii) $S$ is a very large local left order in $Q$.
\end{Proposition}
\begin{proof} The equivalence of (i) and (ii) is inherent in 
Proposition 3.4 of \cite{rings}. For completeness we include a proof here.

\underline{(i) implies (ii)} Let $w\in Q$. Since $Q$ is regular,
$w\,\ar\, e$ for some $e\in E(Q)$ and as $S$ is very large
in $Q$, there is an element $s\in S\cap H_e$. By Lemma 2.1,
$sw\,\eh\, w$. Now $S$ is a weak left order in $Q$, so that
$sw=a\s b$ for some $a,b\in S$. By the remark following Definition 2.3, 
we may assume that $b\leqar a$ in $Q$ so that $b\,\el\, a\s b$ in $Q$
by Lemma 2.1.  Again using the fact that $Q$ is regular,
$b\,\ar\, f$ for some $f\in E(Q)$, and again from the fact
that $S$ is very large in $Q$, we can choose $t\in R_w\cap L_f\cap S$.

In the following egg-box diagram of the $\dee$-class of $w$, a
$*$ in a cell indicates that it is a group $\eh$-class, a convention
we follow elsewhere in this section.

\[\begin{array}{|c|cc|c|c|}
\hline
\begin{array}{cc}
tas&*\\
e&\\
s&\end{array}&&&\begin{array}{cc}
w&tb\\
sw&\end{array}&t\\
\hline
&&&&\\
\hline
as&&&b&
\begin{array}{c}*\\f\end{array}\\
\hline
\end{array}\]
Notice that $tb\in H_w$ by a result  of Green,
(Proposition 2.3.7 of \cite{h}),
and
\[b=asw\,\ar\, as\,\el\, s\]
since $s\,\ar\, sw\leqar a$, so that $tb$ and $as$ are in the
$\eh$-classes as shown.  Again by Proposition 2.3.7 of \cite{h}, $tas\in H_e$ so
that from $tb=tasw$ we have
$w=(tas)\s tb$ where $tas, tb\in S$ and $tas\,\ar\, tb$ in $Q$. Thus (ii)
holds.

\underline{(ii) implies (iii)} For any $q\in Q$, $q=a\s b$
for some $a,b\in S$ with $a\,\ar\, b$ in $Q$. By Lemma 2.1 (v)  we
have that $b\in H_q\cap S$. Thus $S$ is very large in $Q$.

Let $H=H_e$ where $e\in E(Q)$. Since $S$ is very large,
there exists $s\in S\cap H$. Let $q\in H$. Then
$sq=c\s d$ for some $c,d\in S$
with $c\,\ar\, d$ in $Q$,
so that $d\in H$.  Repeated application of Lemma 2.1 (v)
gives
\[scs,s\s c\s s\s \in H_s=H\]
and
\[(scs)(s\s c\s s\s )=ss\s\]
so that $(scs)\s =s\s c\s s\s $. Now
\[q=s\s c\s d=s\s c\s s\s sd=(scs)\s sd\]
and $scs,sd\in H$. Thus $S$ is a local left order in $Q$.

\underline{(iii) implies (i)} Let $q\in Q$.  Let $q'$ be an inverse
of $q$ (in the usual sense of semigroup theory)
and let $s\in H_q\cap S$. By (2) of Theorem 2.2.4 of
[H], there is an inverse $s'$ of $s$ in $H_{q'}$ such that
\[ss'=qq'\mbox{ and }s's=q'q.\]
\[\begin{array}{|c|cc|c|}
\hline
&&&*\\
s&&&qq'\\
q&&&u\,\,\ v\\
\hline
&&&\\
\hline
*&&&s'\\
q'q&&&q'\\
\hline\end{array}\]
By Proposition 2.3.7 of \cite{h}, 
$qs'\in H_{qq'}$, so that
as $S$ is a local left order in $Q$,
\[qs'=u\s v\]
for some $u,v\in H_{qq'}$.
But then
\[q=qq'q=qs's=u\s vs\]
so that $S$ is a weak left order in $Q$.
\end{proof}
\begin{Corollary} Let $Q$ be a regular semigroup with abelian subgroups.
Let $S$ be a subsemigroup of $Q$. Then $S$ is a weak straight left
order in $Q$ if and only if $S$ is a weak straight right order in $Q$.
\end{Corollary}
\begin{proof} Let $S$ be a weak straight left order in $Q$.
By (ii) implies (iii) of
Proposition 3.1, $S$ is a very large local left order
in $Q$. Thus $S\cap H$ is a left order in $H$ for every group
$\eh$-class $H$ of $Q$.  But $H$ is abelian, so that $S\cap H$
is a right order in $H$. This gives that $S$ is a very large
local right order in $Q$, so that the dual of (iii)
implies (ii) of Proposition 3.1 gives that $S$ is
a
weak straight right order in $Q$. 
The converse is dual.
\end{proof}
\begin{Corollary} \cite{vienna} Let $Q$ be a regular semigroup on which $\eh$ is
a congruence. Then any weak left order in $Q$ is straight.
\end{Corollary}
\begin{proof}   By Proposition 3.1 it is enough to show that
any weak left order in $Q$ is very large. 

Let $S$ be a weak left order in $Q$ and let $q\in Q$. Then
$q=a\s b$ for some $a,b\in S$ and as $\eh$ is a congruence
on $Q$,
\[q=a\s b\,\eh\, ab\in S,\]
so that $S$ is very large in $Q$.
\end{proof}

In view of the proliferation of examples of straight left orders
(see \cite{cancell,rings,hii}) and the characterisation of straight
left orders given in the next section, it is natural to ask for
conditions
which ensure straightness, or otherwise, of left orders.  
From the remarks following Definition 2.3, if $S$ is a (two-sided) order
in $Q$, then $S$ intersects every $\ar$-class and every $\el$-class
of $Q$. This latter condition would at first sight
(but erroneously) seem to approach $S$ being
very large in $Q$ and hence, by Proposition 3.1, being straight
in $Q$. Indeed we did not know previously of any orders that
were {\em not} straight.  Here is one such.

\begin{Example}$\, $
\end{Example}  Let $G=\langle a\rangle$ 
be the free cyclic group on $\{ a\}$ and let
$B$ be the Brandt semigroup
 over the trivial group on $\mathbb{Z}\times
\mathbb{Z}$. That is,
\[(u,v)(x,y)=\left\{
\begin{array}{cl}(u,y)&\mbox{ if }v=x\\
0&\mbox{ if }v\neq x.\end{array}\right.\]
Let
\[Q=G\cup B\] and extend the multiplication in $G$ and $B$ by setting
\[a^i(u,v)=(i+u,v)\,\,\, (u,v)a^i=(u,v-i)\]
and
\[a^i0=0=0a^i\]
for all $a^i\in G$ and $(u,v)\in B$.  It is straightforward
to check that $Q$ is a semigroup. Clearly $Q$ is inverse 
with ideal $B$.

Let
\[H=\{ a^i:i\geq 0\}\]
and let
\[C=(\mathbb{N}\times\mathbb{Z})\cup (\mathbb{Z}\times\mathbb{N}^-)\]
where 
\[\mathbb{N}=\{ z\in \mathbb{Z}:z\geq 0\}\mbox{ and }
\mathbb{N}^-=\{ z\in \mathbb{Z}:z<0\} .\]
It is easy to see that no product of elements of $C$ can yield an
element of $\mathbb{N}^-\times\mathbb{N}$, so that $C$ is a subsemigroup
of $B$. It follows that
\[S=H\cup C\]
is a subsemigroup of $Q$. We claim that $S$ is an order in $Q$.

Notice that $Q$ is a monoid with identity $a^0$ and as
$a^0\in S$, for any $q\in S$ we have
\[q=(a^0)\s q=q(a^0)\s.\]
For $i<0$ we have $a^i=(a^{-i})\s a^0$. Consider 
$(u,v)\in B\setminus C$. We have $u<0$ and $v\geq 0$.
Now
\[(u,v)=a^u(0,v)=(a^{-u})\s (0,v)\]
where $a^{-u},(0,v)\in S$ and 
\[(u,v)=(u,-1)a^{-v-1}=(u,-1)(a^{v+1})\s\]
where $a^{v+1},(u,-1)\in S$. Thus $S$ is a weak order in $Q$.

Clearly all the elements of $H$ lie in the subgroup $G$ of $Q$; if
$(u,v)\in S$ is square cancellable, then since
$(u,v)\,\ehs\, (u,v)^2$ in particular $(u,v)^2\neq 0$
so that $u=v$ and $(u,v)\in E(Q)$. Thus $S$ is an order in $Q$.

Since $B$ is an ideal and $\eh$ is trivial on $B$, it follows
that the $\eh$-classes of $Q$ are $G$ and the singletons of $B$.
Thus $S$ is not very large in $Q$ so that by Proposition 3.1
$S$ is not straight in $Q$.

\medskip

We now attempt to isolate the trait of $Q$ in
Example 3.4 that allows it to possess
a non-straight left order.  The crucial  behaviour is the action of
the elements of $G$ on those of $B$.  

Let $Q$ be a semigroup and let $a,q\in Q$ with $a\in\eh(Q)$ and $q\leqar
a$.
Let
\[\rho(a,q)=\{ R_b:b\,\dee\, q,b\leqar a\}.\]
Notice that $R_q\in \rho(a,q)$ and for any $c\in H_a$,
$\rho(a,q)=\rho(c,q)$.  For any $R_b\in\rho(a,q)$,
\[q\,\dee\, b\,\el\, ab\]
by the dual of Lemma 2.1 (i), and clearly $ab\leqar a$, so that
$R_{ab}\in\rho(a,q)$.

\begin{Lemma} With $a,q$ as above,
\[\phi_a:\rho(a,q)\rightarrow\rho(a,q)\]
given by
\[\phi_a(R_b)=R_{ab}\]
is an order automorphism with inverse $\phi_{a\s}$.
\end{Lemma}
\begin{proof} Since $\leqar$ is left compatible,
$\phi_a$ is well defined and order preserving.  Notice
that $\rho(a,q)=\rho(a\s ,q)$ and, for any
$R_b\in\rho(a,q)$,
\[b=aa\s b=a\s ab,\]
whence $\phi_a$ and $\phi_{a\s}$ are mutually inverse order
isomorphisms.
\end{proof}

Putting $\mathcal{A}(X)$ to be the group of order automorphisms
of a partially ordered set $X$, the above shows that
$\phi_a\in\mathcal{A}(\rho(a,q))$. Let $\langle a\rangle$
be the cyclic subgroup of $H_a$ generated by $a$.

\begin{Corollary} With $a,q$ as above,
\[\overline{\phi}=\overline{\phi_{a,q}}:\langle a\rangle
\rightarrow\mathcal{A}(\rho(a,q))\]
given by
\[\overline{\phi}(a^i)=\phi_{a^i}\]
is a morphism.
\end{Corollary}
\begin{proof} As remarked earlier,
\[\rho(a^i,q)=\rho(a,q)=\rho(a^j,q)\]
for any $a^i,a^j\in\langle a\rangle$. Clearly then
\[\phi_{a^i}\phi_{a^j}=\phi_{a^{i+j}},\]
whence $\overline{\phi}$ is a morphism.\end{proof}

 We remark that in the previous result, for any $i\in\mathbb{Z}$
we have
\[\phi_{a^i}=\overline{\phi}(a^i)=(\overline{\phi}(a))^i=\phi_a^i,\]
since $\overline{\phi}$ is a morphism.

Looking again at Example 3.4 we see that
$\rho(a,(0,0))$ consists of all the $\ar$-classes of $B$. Since
\[\phi_{a^i}(0,0)=\phi_{a^j}(0,0) \mbox{ if and only if }i=j,\]
we see that $\overline{\phi}(\langle a\rangle)$
 has infinite order. On the other hand we have the following.

\begin{Lemma} Let $Q$ be a semigroup and let $q\in Q$ be such
that $q=a\s b$ for some $a,b\in Q$ with $b\leqar a$. If
$\overline{\phi}(\langle a\rangle)$ is finite,
where $\overline{\phi}=\overline{\phi_{a,q}}$, then
$q\,\eh a^kb$ for some $k\geq 0$.
\end{Lemma}
\begin{proof} If $\overline{\phi}(\langle a\rangle)$ is finite,
then $\overline{\phi}(a)=\phi_a$ has finite order, say
\[(\phi_a)^n=\phi_{a^n}=I\]
where $n\geq 1$ and $I$ is the identity of $\mathcal{A}(\rho(a,q))$.
Then
\[R_q=I(R_q)=\phi_{a^n}(R_q)=R_{a^nq}\]
so that
\[a\s b\,\ar\, a^na\s b=a^{n-1}b.\]
Since $b\leqar a$ we have
\[q\,\el\, b\,\el\, a^{n-1}b\]
so that $q\,\eh\, a^{n-1}b$ as required.
\end{proof}

The next result follows immediately from  Lemma 3.7.

\begin{Proposition}  Let $S$ be a weak left order in $Q$.
Let $q\in Q$ be such that $q=a\s b$ where
$a,b\in S$ and $b\leqar a$ in $Q$.
If
$\overline{\phi}(\langle a\rangle)$ is finite,
where $\overline{\phi}=\overline{\phi_{a,q}}$,
then $S\cap H_q\neq\emptyset$.
\end{Proposition}

\begin{Corollary} Let $S$ be a weak left order in $Q$. If
for any $q\in Q$ we can write $q=a\s b$ where $a,b\in S$,
$b\leqar a$ in $Q$ and $\rho(a,q)$ is finite, then $S$ is
a straight weak left order in $Q$.
\end{Corollary}
\begin{proof} From Proposition 3.8, $S$ is very large
in $Q$. The result follows by Proposition 3.1.
\end{proof}

\begin{Corollary} Let $S$ be a weak left order in $Q$ such that $Q$ has
finitely many $\ar$-classes. Then $S$ is a straight weak left order in
$Q$.
\end{Corollary}

The conditions of Corollary 3.10 can hold without $\eh$ being
a congruence on $Q$, as the following example shows.

\begin{Example}$\, $
\end{Example} This is a variation on Example 3.4. As in that example,
let $G$ be the free cyclic group on $a$. Let
$n$ be a natural number greater than $2$ and let $B$ be the
Brandt semigroup on $\mathbb{Z}_n\times\mathbb{Z}_n$ with trivial
subgroup. Let $Q=G\cup B$ and extend the multiplication in $G$ and
$B$ by putting
\[a^i(\overline{u},\overline{v})=(\overline{i+u},\overline{v}),
(\overline{u},\overline{v})a^i=(\overline{u},\overline{v-i})\]
and
\[a^i0=0=0a^i\]
for all $a^i\in G$ and $(\overline{u},\overline{v})\in B$. As in
Example 3.4, it is easy to see that $Q$ is an inverse semigroup
with ideal $B$ and finitely many $\ar$-classes.

Let $H=\{ a^i:i\geq 0\}$ and let
$S=H\cup B$. Clearly $S$ is  a subsemigroup of $Q$. For any $q\in B$
we have
\[q=(qq')\s q=q(q'q)\s\]
where $q'$ is the inverse of $q$, and clearly
\[q\,\ar\, qq'\mbox{ and }q\,\el\, q'q.\]
Hence $S$ is a straight weak order in $Q$.
The observation that $\eh$ and $\ehs$ coincide in a regular
semigroup yields easily that $S$ is an order in $Q$.

The relation $\eh$ is not a congruence on $Q$ since $a\,\eh a^0$
but
\[a(\overline{0},\overline{0})=(\overline{1},\overline{0})
\not\hspace{-1mm}\eh (\overline{0},\overline{0})=a^0(\overline{0},\overline{0}).\]

Another situation where Proposition 3.8 applies is as follows.

\begin{Corollary} Let $S$ be a weak left order in a regular
semigroup $Q$ such that in any $\dee$-class of $Q$
the $\ar$-classes form an inversely well ordered chain. Then
$S$ is straight in $Q$.
\end{Corollary}
\begin{proof} Let $q\in Q$. Since $S$ is a weak left order in $Q$,
$q=a\s b$ for some $a,b\in S$ with $b\leqar a$ in $Q$.
Now $\rho(a,q)$ is a non-empty subset of an
inversely well ordered chain, hence is itself an inversely well
ordered chain. It therefore has trivial group of order
automorphisms, so that $S\cap Hq\neq\emptyset$ by Proposition 3.8.
Thus $S$ is very large in $Q$ and the result follows
by Proposition 3.1.
\end{proof}

\begin{Corollary}\cite{hii} Let $S$ be a weak left order in an inverse
$\omega$-semigroup. Then $S$ is straight in $Q$.
\end{Corollary}

Munn shows in \cite{munn} that $\eh$ is a congruence on any inverse
$\omega$-semigroup. Thus the above result is also a consequence of 
Corollary 3.3. 

Finally in this section we give an example of a non-straight
order in an inverse semigroup which has chain of idempotents the
ordinal sum of the dual of $\omega$ with $\omega$.

\begin{Example}$\,$
\end{Example} Let $D=\mathbb{Z}\times\mathbb{Z}$ with
the `bicyclic' multiplication
\[(u,v)(x,y)=(u-v+t,y-x+t)\]
where $t=$ max $\{ v,x\}$. Straightforward calculation
gives that $D$ is a semigroup with idempotents
\[E(D)=\{ (a,a):a\in\mathbb{Z}\}.\]
The idempotents of $D$ form a chain which is the ordinal sum
of the dual of $\omega$ with $\omega$. Moreover, $D$
is regular, hence inverse, where
\[(u,v)'=(v,u)\mbox{ for all }(u,v)\in D.\]

Let $G=\langle a\rangle$ be the free cyclic group on $a$ and let
$Q=G\cup D$. Extend the multiplication in $G$ and $D$ by putting
\[a^i(u,v)=(i+u,v)\mbox{ and }(u,v)a^i=(u,v-i).\]
As in previous examples, it transpires that $Q$ is a semigroup
which is clearly inverse with ideal $D$.

Let
\[H=\{ a^i:i\geq 0\},\,\,  F=\{ (u,v)\in D:u-v\geq 0\}\]
and put
\[S=H\cup F.\]
For any $(u,v),(x,y)\in F$ we have that
\[(u,v)(x,y)=(u-v+t,y-x+t)\]
where $t=$ max $\{ v,x\}$ and
\[(u-v+t)-(y-x+t)=(u-v)+(x-y)\geq 0,\]
so that $F$ is a subsemigroup of $D$. Since also for any
$(u,v)\in F$ and $i\geq 0$ we have
\[u+i\geq u\geq v\geq v-i,\]
$S$ is a subsemigroup of $Q$.  We claim that $S$ is an order in $Q$.

As before we notice that $S$ contains the identity of $Q$, so that
to show that $S$ is a weak order in $Q$ it only remains to express
the elements of $D\setminus F$ as quotients of elements of $S$.

Let $(u,v)\in D\setminus F$, so that $u<v$. Then
\[(u,v)=(v+(u-v),v)=a^{u-v}(v,v)=(a^{v-u})\s (v,v)\]
and $a^{v-u},(v,v)\in S$.  Further,
\[(u,v)=(u,u-(u-v))=(u,u)a^{u-v}=(u,u)(a^{v-u})\s\]
and again, $(u,u),a^{v-u}\in S$. Thus $S$ is a weak 
order in $Q$.

To show that $S$ is an order, it remains to show that 
elements of $F$ that are square cancellable in $S$ lie in
subgroups of $Q$. In fact they are idempotent. For if
$(u,v)\in F\cap\sqs$, then as $u\geq v$ we have
\[(u,v)^2=(2u-v,v)\]
and $2u-v\geq u$. It follows that
\[(u,u)(2u-v,v)=(2u-v,v)=(2u-v,2u-v)(2u-v,v).\]
But $(u,v)\ehs (2u-v,v)$ and so
\[(u,u)(u,v)=(2u-v,2u-v)(u,v).\]
This yields
\[(u,v)=(2u-v,u)\]
whence $u=v$ and $(u,v)\in E(Q)$.

As in the bicyclic semigroup, it is easy to see that $\eh$
is trivial on $D$. Since $D$ is an ideal of $Q$ it follows
that the $\eh$-classes of $Q$ are $G$ and the
singletons of $D$. Thus $S$ fails to be very large in $Q$, hence fails
to be straight.  As in Example 3.4, 
$\rho(a,(0,0))$ consists of all the $\ar$-classes of $D$
and 
the order automorphism $\phi_a$ of $\rho(a,(0,0))$ has
infinite order.

\section{A new characterisation of straight left orders}

The main result of \cite{straight} is a complete characterisation of
those semigroups that are straight left orders. The aim of this
section is to improve upon this result.

The approach of \cite{straight} is via embeddable *-pairs, which
are pairs of preorders whose properties reflect those of
$\leqel$ and $\leqar$ on a regular semigroup and which for convenience
and reference we now detail.

An ordered pair $\pee=(\leq_l,\leq_r)$ of preorders on a semigroup $S$ is a
$*$-{\em pair} if $\leq_l$ is right compatible with
multiplication, $\leq_r$ is left compatible, $\leq_l\,\subseteq\,\leqels$
and $\leq_r\,\subseteq\,\leqars$.  Clearly 
$\pee^*=(\leqels,\leqars)$ is
a $*$-pair for any semigroup $S$. From Lemma 2.2 we have the following.

\begin{Corollary}
Let $S$ be a subsemigroup of $Q$. Then
\[(\leq_{\mathcal{L}^Q}\cap \,(S\times S),
\leq_{\mathcal{R}^Q}\cap \,(S\times S))\]
is a $*$-pair for $S$, the $*$-{\em pair for $S$ induced by $Q$.}
\end{Corollary}

In particular the above applies when $S$ is a left order in $Q$.
Our immediate aim is to describe a number of properties held
by certain induced  $*$-pairs. This will lead us to the notion of
an {\em embeddable $*$-pair}. On the other hand, we know from
\cite{straight} that an embeddable $*$-pair satisfying certain conditions
is {\em always} induced by a semigroup of left quotients. 

For a $*$-pair $\pee=(\leq_l,\leq_r)$, we denote by $\el'$ and $\ar'$ the
equivalence
relations associated with $\leq_l$ and $\leq_r$ respectively. Notice
that $\el'$ is a right congruence contained in $\els$ and $\ar'$
is a left congruence contained in $\ars$; $\eh'$ denotes
the intersection of $\el'$ and $\ar'$ and we put
\[\mathcal{G}(S)=\{ a\in S:a\,\eh'\, a^2\} ,\]
so that $\mathcal{G}(S)\subseteq\sqs$.

In the following lemma, the slightly eccentric notation is that
inherited
from \cite{straight}.

\begin{Lemma}\cite{straight} Let $S$ be very large subsemigroup of a
regular semigroup $Q$. Then the $*$-pair for $S$
induced by $Q$ satisfies the following conditions and the duals
(Eii)(r), (Ev)(r), (Evi)(r) and (Evii)(r)  of
(Eii)(l), (Ev)(l), (Evi)(l) and (Evii)(l) respectively.

(Ei) $\el'\circ\ar'=\ar'\circ\el'$.

(Eii)(l) For all $b,c\in S$, $b\leq_l c$ if and only if $b\,\el'\, dc$
for some $d\in S$.

(Eiii)  Every $\el'$-class and every $\ar'$-class contains an element
from $\mathcal{G}(S)$.

(Ev)(l) For all $a\in\mathcal{G}(S)$ and $b\in S$, if $b\leq_l a$, then
$ba\,\ar'\, b$.

(Evi)(l) For all $a\in\mathcal{G}(S)$ and $b,c\in S$, if $b,c\leq_l a$
and $ba=ca$, then $b=c$.

(Evii)(l) For all $a\in\mathcal{G}(S)$ and $b,c\in S$, if $b,c\leq_l a$
and $ba\,\el'\, ca$, then $b\,\el'\, c$.
\end{Lemma}

As in \cite{straight}  we say that a $*$-pair satisfying conditions
(Ei), (Eii), (Eiii), (Ev), (Evi), and (Evii), where
`(Ex)' means `(Ex)(l) {\em and} (Ex)(r)',  
is an {\em embeddable $*$-pair}.
We make the observation that if 
 a $*$-pair satisfies (Eii)(l) then for any $b,c\in S$,
\[b\leqel c\mbox{ in }S\mbox{ implies that }b\leq_l c\]
and dually, if  (Eii)(r) holds then
\[b\leqar c\mbox{ in }S\mbox{ implies that }b\leq_r c.\]
 If $S$ is a straight left order in $Q$,
then $Q$ is perforce regular and by
Proposition 3.1, $S$ is very large in $Q$; in view of Lemma 4.2
$Q$ induces an embeddable $*$-pair for $S$.
On the other hand, if $S$ is also a straight
left order in $P$, then $Q$ is
isomorphic to $P$ under an isomorphism which restricts to the identity
on $S$, if and only if $Q$ and $P$ induce the {\em same}
embeddable $*$-pair on $S$ \cite{compreg}. Embeddable $*$-pairs
are also crucial in determining the existence of semigroups of
straight left quotients, as the next result shows.

\begin{Theorem}\cite{straight} Let $S$ be a semigroup having an embeddable $*$-pair
$\pee$.
Then
$S$ is a straight left order in a semigroup $Q$ inducing
$\pee$ if and only if $\pee$
satisfies the following conditions.

(Gi) $\sqs =\mathcal{G}(S)$.

(Gii) If $a\in \sqs$ then $H'_a$ is right reversible.

(Giii) If $b,c\in S$ then $b\leq_l c$ if and only if there exist
$h\in \sqs$ and $k\in S$ with $b\leq_r h\,\ar'\, k$ and $hb=kc$.

(Giv) If $a,c\in\sqs$ and $b\in S$ with $a\,\ar'\, b$,
then
there exist $u,h\in \sqs,v,k\in S$
with $h\,\ar'\, k\,\ar'\, u\,\ar'\, v, au\,\ar'\, bc,
v\leq_lc,u\leq_r a,k\leq_l a, hua=ka^2$
and $hvc^2=kbc$.
\end{Theorem}

Of the four conditions listed in the above result, only the first
two appear natural. If $S$ is a straight left order in $Q$ and
$Q$ induces $\pee=(\leq_l,\leq_r)$, then for any $a\in S$,
\[a\in\sqs\mbox{ if and only if }a\,\eh\, a^2\mbox{ in }Q\]
by the definition of left order. Since 
\[\eh^Q\cap (S\times S)=\eh'\]
we have that
\[a\in\sqs\mbox{ if and only if }a\in\mathcal{G}(S),\]
so that (Gi) holds. In view of Proposition 3.1 and 
the theorem of  Ore and Dubreil, 
we certainly have (Gii). Conditions (Giii) and (Giv) are technical
conditions needed in the {\em construction} of semigroups of left
quotients. The object of this section is to show that they hold for
any embeddable $*$-pair satisfying (Gi) and (Gii).
That is, we have the following improved version of Theorem 4.3.

\begin{Theorem} Let $S$ be a semigroup having an embeddable $*$-pair
$\pee$.
Then
$S$ is a straight left order in a semigroup $Q$ inducing
$\pee$ if and only if $\pee$
satisfies conditions (Gi) and (Gii).

\end{Theorem}

\begin{proof} In view of Theorem 4.3 we need only prove the converse.
We begin with some minor but useful observations.

\begin{Lemma} Let $\mathcal{P}$ be a $*$-pair
for a semigroup $S$.  If $\mathcal{P}$ satisfies
(Ev), then for all $a,b\in S$:

(i) if $a\in\mathcal{G}(S)$ and $a\,\el'\, b$, then $ba\,\eh'\, b$;

(ii) if $a\in\mathcal{G}(S)$ and $a\,\ar'\, b$, then $ab\,\eh'\, b$;

(iii) if $a\in\mathcal{G}(S)$, then $H'_a$ is a subsemigroup of
$S$ such that $H'_a\subseteq\mathcal{G}(S)$.
\end{Lemma}
\begin{proof} If $a\in\mathcal{G}(S)$ and $a\,\el'\, b$, then by
(Ev)(l), $ba\,\ar'\, b$. Using the fact that $\el'$ is a right
congruence we also have that
\[ba\,\el'\, a^2\,\eh'\, a\,\el'\, b\]
so that $ba\,\eh'\, b$ and (i) holds. Dually, (ii) holds.

In view of (i) and (ii), if $a\in\mathcal{G}(S)$ and $b,c\in H'_a$,
we have that
\[bc\,\ar'\, ba\,\eh'\, b\,\eh'\, a\]
and
\[bc\,\el'\, ac\,\eh'\, c\,\eh'\, a\]
so that $bc\in H'_a$ and $H_a'$ is a subsemigroup. 
Clearly then $H'_a\subseteq\mathcal{G}(S)$.
\end{proof}

We can now deduce a result strongly remiscent of
a classical result of Green for regular $\dee$-classes.

\begin{Corollary}  Let $\mathcal{P}$ be a
$*$-pair for a semigroup $S$. If $\mathcal{P}$ satisfies (Ev), then
for all $u,v\in S$
\[u\,\el'\, s\,\ar'\, v\mbox{ where }s\in\mathcal{G}(S),\]
implies that
\[u\,\ar'\, uv\,\el'\, v.\]
\end{Corollary}
\begin{proof}  If $u,v,s\in S$ are as given,
then
\[uv\,\el'\, sv\,\eh'\, v\]
by Lemma 4.5. Dually, $uv\,\ar'\, u$.
\end{proof}

The following lemma corresponds to a stronger version of condition
(Giii).

\begin{Lemma} Let $\mathcal{P}$ be a $*$-pair for a
semigroup
$S$. Suppose that $\mathcal{P}$ satisfies conditions
\begin{center}(Ei), (Eii), (Eiii), (Ev), (Evi)(l) and (Gii).\end{center}
Then for any $b,c\in S$,
\[b\leq_l c\mbox{ if and only if }hb=kc\]
for some $h\in \mathcal{G}(S),k\in S$ with
\[h\,\ar'\, k\,\ar'\, b.\]
\end{Lemma}
\begin{proof}  First, if $h$ and $k$ exist as given, then from (Ev)(r),
\[b\,\el'\, hb=kc\leq_l c,\]
the latter relation by (Eii)(l).

Conversely, suppose that $b\leq_l c$. Again by (Eii)(l), 
$b\,\el'\, dc$ for some $d\in S$. By
(Eiii), $R'_b,L'_b$ and $R'_{dc}$ contain elements of $\mathcal{G}(S)$.

In the following diagram, and others elsewhere in this section,
the rows correspond to the $\ar'$-classes, the columns to the
$\el'$-classes, and a $*$ in a cell, that is, in an
$\eh'$-class, indicates that it consists of elements of $\mathcal{G}(S)$.

\[\begin{array}{|c|cc|c|c|}
\hline
\begin{array}{cc}
by&*\\
xdcy&\\
u&v\end{array}
&&&\begin{array}{c}b\\xdc\end{array}&x\\ \hline
&&&&\\ \hline
&&&dc&*\\
\hline
y&&& *&\\
\hline
\end{array}
\]
By (Ei) we can choose $x,y\in S$ as shown. Calling upon Corollary 4.6,
we
have
\[b\,\ar'\, by\,\el'\, y\]
and
\[x\,\ar'\, xdc\, \el'\, dc\]
so that $xdc\,\eh'\, b$ and again by Corollary 4.6,
$xdcy\,\eh'\, by$. By Lemma 4.5, $H'_{by}\subseteq\mathcal{G}(S)$
and $H'_{by}$ is right reversible by (Gii). Thus
\[uby=vxdcy\]
for some $u,v\in H'_{by}$ and so
\[(ub)(yb)=(vxdc)(yb).\]
By Lemma 4.5, $ub,vxdc\in H'_b$ and by
Corollary 4.6, $yb\,\el'\, b$. Since $yb\in\mathcal{G}(S)$
(Evi)(l) gives
that 
\[ub=vxdc.\]
Using (Eii)(r) we have that
\[vxd\leq_r v\,\ar'\, ub=vxdc\leq_r vxd\]
so that
\[u\,\ar'\, b\,\ar'\, vxd.\]
Putting $h=u$ and $k=vxd$, the lemma follows.
\end{proof}

It is worthwhile mentioning that in the above argument,
we can arrange $h$ to be $\eh'$-related to {\em any}
 element of $\mathcal{G}(S)$ in the $\ar'$-class of $b$.

\begin{Corollary}
 Let $\mathcal{P}$ be a $*$-pair for a
semigroup
$S$. Suppose that $\mathcal{P}$ satisfies conditions
\begin{center}(Ei), (Eii), (Eiii), (Ev), (Evi)(l) and (Gii).\end{center}
Then for any $b,c\in S$,
\[b\leq_l c\mbox{ if and only if }hb=kc\]
for some $h\in\mathcal{G}(S),k\in S$ with
\[b\leq_r h\,\ar'\, k.\]
\end{Corollary}
\begin{proof} If $h$ and $k$ exist as given, then as in Lemma 4.7,
$b\leq_l c$.  Conversely, if $b\leq_lc$
then according to Lemma 4.7, $h$ and $k$ exist as given, actually
satisfying the stronger condition of being in the $\ar'$-class of
$b$.
\end{proof}

The next lemma is reminiscent of Lemma 2.1(iv).
\begin{Lemma} Let $\mathcal{P}$ be
a $*$-pair for a semigroup $S$ such that (Eii) holds. If
$p\in S$ and $q\in \mathcal{G}(S)$ then
\[q\leq_r p\mbox{ implies that }qp\,\ar'\, q\]
and dually,
\[q\leq_l p\mbox{ implies that }pq\,\el'\, q.\]
\end{Lemma} 
\begin{proof} If $q\leq_r p$, then since $\leq_r$ is left compatible,
\[q\,\eh'\, q^2\leq_r qp\leq_r q,\]
the latter relation by (Eii)(r). Thus $q\,\ar'\, qp$.
\end{proof}

Returning to the proof of Theorem 4.4, suppose that
$\mathcal{P}=(\leq_l,
\leq_r)$ is an embeddable $*$-pair for $S$ satisfying (Gi) and (Gii).
Condition
(Giii) is now an immediate consequence of (Gi) and Corollary 4.8. It
remains to show that (Giv) holds.

Let $a,c\in \sqs=\mathcal{G}(S)$ and let $b\in S$ be such that
$a\,\ar'\, b$. We have by (Eii),
\[bc\leq_l c\,\eh'\, c^3\]
so that by Lemma 4.7,
\[pbc=qc^3\]
for some $p,q\in S$ with $bc\,\ar'\, p\,\ar'\, q$ and
$p\in\sqs$. Since $p\,\ar'\, bc\leq_r a$ and $p\in\sqs$,
$p\,\ar' pa$ by Lemma 4.9.
\[
\begin{array}{|c|c|cc|c|c|}
\hline
\begin{array}{c}
bc\\
pbc=qc^3\end{array}
&a^2t&&&\begin{array}{cc}
p&*\\
as&\end{array}&pa\\
\hline
&&&&&\\
\hline
&at&&&s&spa\\
\hline
\end{array}\]
We have
\[bc\leq_rb\,\ar'\, a^2\]
so that by (Eii),
\[bc\,\ar'\, a^2t\]
for some $t\in S$. Now $at\leq_r a$ so that by (Ev), $a^2t\,\el'\, at$.
By (Ei) we can choose $s$ with $at\,\ar'\, s\, \el'\, p$. By Corollary
4.6,
\[s\,\ar'\, spa\,\el'\, pa.\]
We claim that $spa\in\mathcal{G}(S)$.

First, as $s\,\ar'\, at\leq_r a$ we have using (Ev)(r) that
\[p\,\el'\, s\, \el'\, as.\]  But $as\,\ar'\, a^2t\,\ar'\, p$ so that
$as\,\eh'\, p$.  Now
\[(spa)^2=(spa)(spa)=(spas)pa\,\ar'\, (spas)p=s(pasp).\]
But $as\,\eh'\, p$ and $H'_p$ is a subsemigroup by Lemma 4.5,
so that
\[(spa)^2\,\ar'\, sp\,\ar'\, spa.\]
On the other hand,
\[(spa)^2=s(paspa)\,\el'\, p^2aspa\,\el'\, pa\,\el'\, spa.\]
Hence $spa\,\eh'\, (spa)^2$ and $spa\in\mathcal{G}(S)$ as claimed.

Let $u=spa$. By (Gi) $u\in\sqs$ and
\[au=aspa\,\ar'\, a^2t\,\ar'\, bc.\]
Further, as $u\,\ar'\, at\leq_r a$, Lemma 4.9 gives that
\[ua\,\ar'\, u.\]
But $ua\leq_l a\,\el'\, a^3$ so that by Lemma 4.7,
\[hua=k'a^3\]
for some $h\in\sqs$ and $k'\in S$  with
\[h\,\ar'\, k'\,\ar'\, ua\, \ar'\, u.\]
Put $k=k'a$. Then $k\leq_l a$ and
\[hua=ka^2.\]
Moreover,
\[k=k'a\,\ar'\, k'a^3=hua\,\eh'\, ua\,\ar' h.\]
Replacing $u$ by $spa$ we have
\[hspa^2=hua=ka^2\]
so that as $a^2\,\ar'\, b$ and $\ar'\subseteq \ars$,
\[hspb=kb.\]
Recall that $pbc=qc^3$ so that
\[hsqc^3=hspbc=kbc.\]
Put $v=sqc$.  Then $v\leq _lc$ and
\[hvc^2=kbc.\]
Finally,
\[v=sqc\,\ar'\, sqc^3=spbc.\]
But
\[pbc\,\eh'\, bc\,\ar'\, p\]
so that
\[v\,\ar'\, spbc\,\ar'\, sp\,\eh'\, s\,\ar'\, spa=u.\]
This completes the proof of Theorem 4.4.
\end{proof}

We can immediately deduce the result in the two-sided case.

\begin{Corollary} Let $S$ be a semigroup having embeddable $*$-pair
$\pee$. Then $S$ is a straight  order in a semigroup
$Q$ inducing $\pee$ if and only if $\pee$ satisfies
conditions (Gi), (Gii) and the dual (Gii)' of (Gii).
\end{Corollary}
\begin{proof} The direct implication is immediate from Theorem 4.4
and its dual.
Conversely, if (Gi), (Gii) and (Gii)' hold, then the same result
gives that $S$ is a straight left order in a semigroup $Q$ inducing
$(\leq_l,\leq_r)$.

Let $q\in H$ where $H=H_q$ is a group $\eh$-class of $Q$. By
Proposition 3.1, $S$ is a local left order in $Q$ so that $q=a\s b$
for some $a,b\in S\cap H$.  Since $Q$ induces $(\leq_l,\leq_r)$,
$S\cap H=H'_a$ and $H'_a$ is left reversible by (Gii)'. Thus
$ac=bd$ for some $c,d\in S\cap H$ so that $cd\s =a\s b=q$
and $S$ is a local right order in $Q$. As $S$ is very large
in $Q$, $S$ is a straight right order in $Q$ by the dual
of Proposition 3.1.
\end{proof}

With some care we can improve upon the previous result. Condition (Evii)
in the definition of embeddable $*$-pair becomes superfluous.

\begin{Lemma} Let $\mathcal{P}$ be a $*$-pair
for a semigroup $S$. If $\mathcal{P}$ satisfies
\begin{center} (Ei),(Eii), (Eiii), (Ev), (Evi)(l) and (Gii),\end{center}
then $\mathcal{P}$ also satisfies (Evii)(l).
\end{Lemma}
\begin{proof} Let $b,c\in S$ and $a\in\mathcal{G}(S)$ with
$b,c\leq_l a$ and $ba\,\el'\, ca$. By Lemma 4.7,
\[hba=kca\]
for some $h\in\mathcal{G}(S)$ and $k\in S$ with $h\,\ar'\, k\,\ar'\,
ba$.
By (Eii)(l),
\[hb\leq_l b\leq_la\mbox{ and }kc\leq_l c\leq_l a\]
so that by (Evi)(l),
\[hb=kc.\]
By (Ev)(l), $ba\,\ar'\, b$ so that $b\,\ar'\, h$ and by (Ev)(r),
\[b\,\el'\, hb=kc\leq_l c,\]
the latter relation by (Eii)(l). Thus $b\leq_l c$
and dually, $c\leq_l b$, so that $b\,\el'\, c$ as required.
\end{proof}

\section{Specialisations}

Theorem 4.4 can be used in a number of ways. As in \cite{straight} we
can recover the
description of straight left orders in regular semigroups
 on which Green's relation $\eh$ is a
congruence.
On the other hand we could use Theorem 4.4 to obtain known descriptions
of straight left orders in more specific classes of semigroups, for
example,
the class of inverse  $\omega$-semigroups (for which indeed $\eh$ is a
congruence \cite{munn}).  The arguments in these cases would of course
be now more succinct, as there is no need to verify conditions (Giii) and
(Giv).

Just as large classes of left orders are naturally straight, so are many
naturally {\em fully stratified}, where a left order $S$ in $Q$
is  fully stratified if $Q$ induces the $*$-pair $\mathcal{P}^*=(\leqels,\leqars)$.
Theorem 4.4 can be specialised to describe fully stratified straight
left orders. In this case, not surprisingly, we can streamline the
conditions on $\mathcal{P}^*$ for it to be embeddable. In the case where
every $\els$-class and every $\ars$-class contains an idempotent, that
is, $S$ is {\em abundant} \cite{f}, we can make even more progress. These results
are contained in a sequel \cite{appl}.

A classic tool in structure theory of semigroups is to decompose
them into `manageable' parts and then `glue'
the parts back together. We took this line in \cite{compreg}, where we
describe straight left orders in completely regular
semigroups in terms of straight left orders in completely simple
semigroups.  With our improved description of straight left orders
we can pursue this approach for wider classes than completely regular
semigroups. We return briefly to this topic at
the end of this section. The details are in \cite{poorders}.

In this paper we content ourselves by giving three straightforward
corollaries of Theorem 4.4.

\begin{Corollary} Let $S$ be a semigroup having an embeddable $*$-pair
$\pee$. Then $S$ is a straight left order
in an inverse semigroup $Q$ inducing $\mathcal{P}$ if and only if
$\pee$ satisfies conditions (Gi), (Gii) and (I):

(I): for all $a,h\in\mathcal{G}(S)$,
\[a\,\ar'\, h\mbox{ implies that }a\,\eh'\, h\]
and
\[a\,\el'\, h\mbox{ implies that }a\,\eh'\, h.\]
\end{Corollary}
\begin{proof} If $S$ is a straight left order in an inverse
semigroup $Q$ inducing $\mathcal{P}$, then certainly (Gi) and
(Gii)
hold by Theorem 4.4. If $a,h$ are $\ar'$-related elements
of $\mathcal{G}(S)$, then $a,h\in \sqs$ and hence $a,h\in\eh(Q)$ since
$S$ is a left order in $Q$. But $Q$ induces $\mathcal{P}$ so that
\[aa\s\,\eh\, a\,\ar\, h\,\eh\, hh\s\]
in $Q$. But $Q$ is inverse so that
\[a\,\eh\,aa\s=hh\s\,\eh\, h\]
in $Q$, giving that $a\,\eh'\, h$ in $S$. Together with its dual this
argument gives that (I) holds.

Conversely, suppose that $S$ satisfies (Gi), (Gii) and (I), so
that by Theorem 4.4, $S$ is a straight left order in $Q$ (perforce
regular) where $Q$ induces $\mathcal{P}$.  Let $e,f\in E(Q)$ be
$\ar$-related in $Q$. Since $S$ is very large in $Q$, $e\,\eh\, a$ and
$f\,\eh\, h$ for some
$a,h\in\eh(Q)\cap S=\mathcal{G}(S)$. Then
\[a\,\eh\, e\,\ar\, f\,\eh\, h\]
so that $a\,\ar'\, h$ in $S$ and by (I), $a\,\eh'\, h$.
Hence 
\[e\,\eh\, a\,\eh\, h\,\eh\, f\]
so that $e=f$.  Thus every $\ar$-class of $Q$ contains a 
unique idempotent. Since the dual result holds for
$\el$-classes, $Q$ is inverse.
\end{proof}

\begin{Corollary}  Let $S$ be a semigroup having an embeddable $*$-pair
$\mathcal{P}$. Then $S$ is a straight
left order in a completely regular semigroup $Q$ such that $Q$ induces
$\mathcal{P}$ if and only if $\pee$ satisfies conditions (GI) and (Gii),
where:

(GI) $\mathcal{G}(S)=S$.
\end{Corollary} 
\begin{proof} If $S$ is a straight left order in a completely regular
semigroup $Q$ inducing
$\pee$, then (Gi) and (Gii) hold by Theorem 4.4.
But $\sqs=S$ since $Q$ is a union of groups, so that
$\mathcal{G}(S)=S$ and (GI) holds.

Conversely, if (GI) and (Gii) hold, then as
$\eh'\subseteq\ehs$,
\[S=\mathcal{G}(S)\subseteq\sqs\subseteq S,\]
so that $\mathcal{G}(S)=\sqs$ and (Gi) holds.
Thus $S$ is a straight left order in a (regular) semigroup
$Q$ inducing $\pee$.   For any $q\in Q$,
$q\,\eh\, a$ for some $a\in S$, since $S$ is very large
in $Q$. But $a\in\sqs$ so that $a$, and hence $q$, lie in
a subgroup of $Q$, which is therefore completely regular.
\end{proof}

We remark that in the previous corollary, condition (Eiii) can be
omitted in the definition of embeddable $*$-pair, since this
is a consequence of (GI).

For our final results we introduce relations $\dee',\leq_j$ and $\jay'$.
These are clearly meant to correspond to $\dee,\leq_{\mathcal{J}}$ and
$\jay$, where for elements $p,q$ of a regular semigroup $Q$,
\[p\leqjay q\mbox{ if and only if }p=xqy\mbox{ for some }x,y\in Q.\]

Let $S$ be a semigroup having embeddable $*$-pair
$\pee=(\leq_l,\leq_r)$.
The relations $\dee',\leq_j$ and $\jay'$ are defined
on $S$ by 
\[\dee'=\el'\vee\ar',\]
for all $a,b\in S$,
\[a\leq_j b\mbox{ if and only if }a\dee' xby\]
and 
\[a\,\jay'\, b\mbox{ if and only if }a\leq_j b\leq_j a.\]

\begin{Lemma} Let $\pee$ be a $*$-pair for a
semigroup $S$.  If $\pee$ satisfies
\begin{center}(Ei),(Eii),(Eiii) and (Ev),\end{center}
then

(i) $\dee'=\el'\circ\ar'$;

(ii) for any $b\in S$,
\[J'(b)=\{ s\in S:s\leq_j b\}\]
is an ideal of $S$ containing $b$;

(iii) $\leq_j$ is a preorder on $S$.
\end{Lemma}
\begin{proof} The first statement is immediate from (Ei).

We next remark that $\leq_j$ is reflexive. For if $b\in S$, then
by (Eiii), we can find $x,y\in \mathcal{G}(S)$ with
$x\,\ar'\, b\,\el'\, y$. By Lemma 4.5,
$b\,\eh'\, xby$ so that certainly $b\,\dee'\, xby$, $b\leq_j b$ and
$b\in J'(b)$.

Again, let $b\in S$ and let $s\in J'(b)$ and $t\in S$. Since $s\in
J'(b)$,
$s\leq_j b$ so that $s\,\dee'\, ubv$ for some $u,v\in S$. By (i),
\[s\,\ar'\, w\,\el'\, ubv\]
for some $w\in S$ so that as $\ar'$ is a left congruence,
\[ts\,\ar'\, tw\leq_lw\,\el'\, ubv,\]
the inequality by (Eii)(l). The same condition gives
that $tw\,\el'\, zubv$ for some $z\in S$ and so
$ts\,\dee'\, zubv$ and
$ts\leq_j b$. That is, $ts\in J'(b)$ so that $J'(b)$ is a left ideal
containing
$b$.   Dually, $J'(b)$ is a right ideal.

It remains only to show that $\leq_j$ is transitive.  Suppose that
$b,c,d\in S$
and
\[b\leq_j c\leq_j d.\]
Then $b\,\dee'\, hck$ for some $h,k\in S$ and $c\in J'(d)$.  But
$J'(d)$ is an ideal, so that $hck\in J'(d)$ and
\[b\,\dee'\, hck\,\dee'\, mdn\]
for some $m,n\in S$. Thus $b\leq_j d$ as required.
\end{proof}

\begin{Corollary}  Let $\pee$ be a $*$-pair for
a semigroup $S$ satisfying
\begin{center} (Ei),(Eii),(Eiii) and (Ev).\end{center}
Then

(i) $\jay'$ is an equivalence on $S$ and the $\jay'$-classes are
partially
ordered by the relation induced by $\leq_j$;

(ii) $\dee'\subseteq\jay'$.
\end{Corollary}
\begin{proof} The first statement is standard. For the second,
suppose that $a,b\in S$ and $a\,\dee'\, b$. As in Lemma 5.3, 
$b\,\dee'\, xby$ for some $x,y\in S$, so that
$a\,\dee'\, xby$ and $a\leq_j b$. Dually, $b\leq_j a$ so that
$a\,\jay'\, b$.
\end{proof}

\begin{Lemma} Let $\pee$ be an embeddable $*$-pair for
a semigroup $S$ and suppose that $S$ is a straight left order in $Q$
inducing
$\pee$. Then

(i) $\dee^Q\cap(S\times S)=\dee'$;

(ii) $\leq_{\mathcal{J}^Q}\cap \,(S\times S)= \, \leq_j$;\\
and

(iii) $\jay^Q\cap (S\times S)=\jay'$.
\end{Lemma}
\begin{proof} The first statement follows immediately from the fact
that $S$ is very large in $Q$.

Suppose now that $a,b\in S$ and $a\leq_j b$, so that
$a\,\dee'\, xby$ for some $x,y\in S$. By (i),
$a\,\dee \, xby$ in $Q$ so that certainly $a\,\jay\, xby$ 
and $a\leq_{\mathcal{J}} b$ in $Q$.

Conversely, suppose that $a,b\in S$ and $a\leq_{\mathcal{J}} b$ in $Q$.
Then $a=pbq$ for some $p,q\in Q$ and as $S$ is very large in $Q$,
$p\,\eh\, c$ and
$q\,\eh\, d$ for some $c,d\in S$. Then
\[a=pbq\,\ar\, pbd\,\el\, cbd\]
so that $a\,\dee\, cbd$ in $Q$. By (i), $a\,\dee'\, cbd$
in $S$ so that $a\leq_j b$. Hence (ii) holds; (iii) follows immediately.
\end{proof}

The proof of the following corollary is immediate from Lemma 5.5
and the fact that $S$ is very large in $Q$.

\begin{Corollary} Let $\pee$ be an embeddable
$*$-pair for a semigroup $S$ and suppose that $S$ is a straight left
order in $Q$ inducing $\pee$. Then

(i) $Q$ is simple if and only if $\jay'=S\times S$;

(ii) $Q$ is bisimple if and only if $\dee'=S\times S$.
\end{Corollary}

Our final aim is to characterise straight left orders in completely
semisimple semigroups. We recall from \cite{cp} that a semigroup $Q$ is
completely
semisimple if and only if its principal factors are either completely
0-simple or completely simple.

In the spirit of the notation of \cite{cp} we define conditions $M^*_r$
and
$M^*_l$ on a $*$-pair $\pee=(\leq_l,\leq_r)$ on a semigroup $S$
satisfying (Ei), (Eii), (Eiii) and (Ev). From Corollary 5.4, each
$\jay'$-class is a union of $\el'$-classes and of
$\ar'$-classes. We say that $\pee$ has $M^*_r$ if, for any
$\jay'$-class, the $\ar'$-classes
it contains satisfy the descending chain condition.  Condition $M^*_l$
is defined dually.  

\begin{Theorem} Let $S$ be a straight left order in a semigroup $Q$
inducing the embeddable $*$-pair $\pee$ for $S$. Then
the
following conditions are equivalent:

(i) $Q$ is completely semisimple;

(ii) for any $a,b\in S$,
\[a\,\jay'\, ab\mbox{ implies that }a\,\ar'\, ab\]
and
\[a\,\jay'\, ba\mbox{ implies that }a\,\el'\, ba;\]

(iii) $\pee$ satisfies $M^*_r$ and $M^*_l$.
\end{Theorem}
\begin{proof} \underline{(i) implies (ii)}  Let $a,b\in S$. If
$a\,\jay'\, ab$, then by Lemma 5.5,
$a\,\jay\, ab$ in $Q$, so that $a$ and $ab$
are non-zero elements of a completely (0-)simple principal factor
$F$ of $Q$. If $ab\,\el\, f\in E(F)$, then
\[ab=(ab)f=a(bf)\]
so that in $Q$,
\[ab\leq_{\mathcal{J}} bf\leq_{\mathcal{J}} f\,\jay\, ab\]
so that $bf\in F\setminus\{ 0\}$.  Thus
$a,a(bf)\in F\setminus\{ 0\}$ and as
$F$ is completely (0-)simple, $a\,\ar\, a(bf)$ in $F$ and hence
$a\,\ar\, a(bf)=ab$ in $Q$.  Thus $a\,\ar'\, ab$ in $S$. Together
with the dual argument, this gives that (ii) holds.

 \underline{(ii) implies (iii)} Let $J'$ be the $\jay'$-class of an
element $a$ of $S$ and suppose that $b\leq_r a$ for some $b\in J'$.
By (Eii)(r), $b\,\ar'\, ac$ for some $c\in S$ and as
$\ar'\subseteq\dee'\subseteq\jay'$,
we certainly have that $a\,\jay'\, ac$ so that by (ii),
\[b\,\ar'\, ac\,\ar'\, a.\]
Thus there are no $\leq_r$ chains of length greater than
one in $J'$ so that certainly $M^*_r$ and dually, $M^*_l$, hold.

 \underline{(iii) implies (i)}  Let $J$ be a $\jay$-class of $Q$. Suppose
that $J$ contained an infinite chain
\[q_1>_{\mathcal{R}}q_2>_{\mathcal{R}}>q_3\hdots \]
Since $S$ is very large in $Q$ we can choose elements $a_i\in H_{q_i},
i\in \mathbb{N}$ and obtain a chain
\[a_1>_{\mathcal{R}}a_2>_{\mathcal{R}}>a_3\hdots \]
so that
\[a_1>_ra_2>_r>a_3\hdots \]
But the $a_i$'s are all $\jay$-related in $Q$, hence $\jay'$-related
in $S$  by Lemma
5.5. This contradicts $M_r^*$. Thus no such chain of $q_i$'s exists.
Consequently,
$J$ contains a minimal $\ar$-class. Dually, $J$ contains a minimal
$\el$-class. Since $Q$ is regular, and so certainly semisimple,
Theorem 6.45 of \cite{cp} gives that $Q$ is completely semisimple.
\end{proof} 

We finish with some remarks which outline an alternative approach
to Theorem 5.7.

Let $S$ be a semigroup having $*$-pair $\pee=(\leq_l,\leq_r)$ satisfying
(Ei),(Eii),(Eiii) and (Ev).  Let $b\in S$. In view of (ii)
and (iii) of Lemma 5.3,
\[I'(b)=\{ s\in S:s<_j b\}\]
is empty or is an ideal of $S$. The $\pee$-{\em principal factors}
of $S$ are the Rees quotients
\[J'(b)/I'(b).\]
It is easy to see that in Theorem 5.7, the $\pee$-principal factors
of $S$ are straight left orders in the principal factors of $Q$, 
which are completely (0-)simple. Fountain and Petrich
 have described left orders in
completely (0-)simple semigroups
(see \cite{fpi} for the two-sided case,
the result for left orders can be deduced from
Theorem  5.7 above). Since $\eh$ is a congruence on
a completely (0-)simple semigroup, such left orders are always
straight. The structure of completely (0-)simple semigroups
is tightly fixed by Rees' theorem.
It is then natural to ask, given a semigroup $S$ which can be
`sliced up' into Rees quotients that are left orders in completely
(0-)simple semigroups, under which conditions is $S$ a left order in a
completely
semisimple semigroup? What if the `slices' are left orders in semigroups
from other well known classes? We use embeddable $*$-pairs and
Theorem 4.4 to investigate these questions in \cite{poorders}.

\end{document}